\documentclass[11pt]{elsarticle}
\usepackage{hyperref}
\usepackage[textwidth=16cm, textheight=24cm, left=2.5cm, right=2.5cm, top=3cm, bottom=3cm]{geometry}
\usepackage{latexsym}
\usepackage{enumerate}
\usepackage{amsmath, amsthm}
\usepackage{amssymb}
\usepackage{graphics}
\usepackage{graphicx}
\usepackage{subfigure}
\usepackage{indentfirst}
\usepackage{mathrsfs}
\usepackage{tikz}
\usepackage{longtable}
\usepackage{multirow}

\usetikzlibrary{decorations.pathreplacing}

\theoremstyle{plain}

\theoremstyle{definition}

\numberwithin{equation}{section}

\makeatletter
\def\ps@pprintTitle{%
  \let\@oddhead\@empty
  \let\@evenhead\@empty
  \def\@oddfoot{\reset@font\hfil\thepage\hfil}
  \let\@evenfoot\@oddfoot
}
\makeatother

\begin{document}
\begin{frontmatter}
\title{Some new trace formulas of tensors with applications in spectral hypergraph theory\tnoteref{title1}}
\tnotetext[title1]{Research supported by the National Science Foundation of China No.11231004, 11271288 and the Hong Kong Research Grant Council (Grant No. PolyU 501909, 502510, 502111 and 501212)}
\author[a]{Jia-Yu Shao}
\ead{jyshao@tongji.edu.cn}
\author[b]{Liqun Qi}
\ead{maqilq@polyu.edu.hk}
\author[b]{Shenglong Hu}
\ead{Tim.Hu@connect.polyu.hk}

\address[a]{Department of Mathematics,  Tongji University,  Shanghai,  China}
\address[b]{Department of Applied Mathematics, The Hong Kong Polytechnic University, Kowloon, Hong Kong}

\date{}
\begin{abstract} We give some graph theoretical formulas for the trace $Tr_k(\mathbb {T})$ of a tensor $\mathbb {T}$ which do not involve the differential operators and auxiliary matrix. As applications of these trace formulas in the study of the spectra of uniform hypergraphs, we give a characterization (in terms of the traces of the adjacency tensors) of the $k$-uniform hypergraphs whose spectra are $k$-symmetric, thus give an answer to a question raised in [3]. We generalize the results in [3, Theorem 4.2] and [5, Proposition 3.1] about the $k$-symmetry of the spectrum of a $k$-uniform hypergraph, and answer a question in [5] about the relation between the Laplacian and signless Laplacian spectra of a $k$-uniform hypergraph when $k$ is odd. We also give a simplified proof of an expression for $Tr_2(\mathbb {T})$ and discuss the expression for $Tr_3(\mathbb {T})$.

\end{abstract}

\end{frontmatter}

\section{Introduction}

\vskip 0.1cm

As was in [9], an order $m$ dimension $n$ tensor $\mathbb {A}=(a_{i_1i_2\cdots i_m})_{1\le i_j\le n \ (j=1,\cdots ,m)}$ over the complex field $\mathbb {C}$ is a multidimensional array with all entries $a_{i_1i_2\cdots i_m}\in \mathbb {C} \  \ (i_1,\cdots ,i_m\in [n]=\{1,\cdots ,n\})$. For a vector $x=(x_1,\cdots ,x_n)^T\in \mathbb {C}^n$, let $\mathbb {A}x^{m-1}$ be a vector in $\mathbb {C}^n$ whose $i$th component is defined as the following:
$$(\mathbb {A}x^{m-1})_i=\sum_{i_2,\cdots ,i_m=1}^n a_{ii_2\cdots i_m}x_{i_2}\cdots x_{i_m} \eqno {(1.1)}$$
and let $x^{[r]}=(x_1^r,\cdots ,x_n^r)^T$. Then ([2, 9]) a number $\lambda \in \mathbb {C}$ is called an eigenvalue of the tensor $\mathbb {A}$ if there exists a nonzero vector $x\in \mathbb {C}^n$ such that
$$\mathbb {A}x^{m-1}=\lambda x^{[m-1]} \eqno {(1.2)}$$
and in this case, $x$ is called an eigenvector of $\mathbb {A}$ corresponding to the eigenvalue $\lambda$. L.Qi and L.H.Lim also defined several other types of eigenvalues (and eigenvectors) in [7, 9].

\vskip 0.28cm

\vskip 0.2cm

The unit tensor of order $m$ and dimension $n$ is the tensor $\mathbb {I}=(\delta _{i_1,i_2,\cdots ,i_m})$ with entries as follows:
$$\delta _{i_1,i_2,\cdots ,i_m}=\left\{
                            \begin{array}{cc}
                              1 &   \mbox{if } i_1=i_2=\cdots =i_m  \\
                              0 &  \mbox{otherwise}\\
                            \end{array}
                          \right.$$
It is easy to see from the definition that $\mathbb {I}x^{m-1}=x^{[m-1]}$. Thus equation (1.2) can be rewritten as
$$(\lambda \mathbb {I}- \mathbb {A})x^{m-1}=0  $$

\vskip 0.1cm

By using the definition of determinants of tensors, Qi ([9]) defined the characteristic polynomial of a tensor $\mathbb {A}$ as the determinant $Det(\lambda \mathbb {I}- \mathbb {A})$, where $\mathbb {I}$ is the unit tensor.

\vskip 0.38cm

In [8], Morozov and Shakirov defined the $d$-th order trace $Tr_d(\mathbb {T})$ of a tensor $\mathbb {T}$ (with order $m$ and dimension $n$), in terms of an auxiliary matrix $A$ of order $n$ (whose entries $a_{ij}$ are viewed as independent variables) and some differential operators involving these variables as follows:

$$Tr_d(\mathbb {T})=(m-1)^{n-1}\sum_{d_1+\cdots + d_n=d}\prod_{i=1}^n \frac {1} {(d_i(m-1))!}\left (\sum_{y_i\in [n]^{m-1}}t_{iy_i}\frac {\partial} {\partial a_{iy_i}}\right )^{d_i}tr(A^{d(m-1)}) \eqno {(1.3)}$$
where we use the notations
$$t_{iy}:=t_{ii_2\cdots i_m} \qquad \mbox {and} \qquad \frac {\partial} {\partial a_{iy}}:=\frac {\partial} {\partial a_{ii_2}}\cdots \frac {\partial} {\partial a_{ii_m}} \qquad \mbox {(if} \ y=i_2\cdots i_m )$$
and $d_1,\cdots , d_n$ run over all nonnegative integers with $d_1+\cdots + d_n=d$.

\vskip 0.38cm

It was proved in [4, Theorem 6.3] that, by using $Tr_d(\mathbb {T})$ and the Schur function, the characteristic polynomial $\phi_{\mathbb {T}}(z)$ of $\mathbb {T}$  can be expressed in the following way:

$$\phi_{\mathbb {T}}(z)= \sum_{k=0}^{d}P_k\left (-\frac {Tr_1(\mathbb {T})} {1},\cdots, -\frac {Tr_k(\mathbb {T})} {k}\right )z^{d-k}  \qquad (d=n(m-1)^{n-1})  $$

\vskip 0.18cm

\noindent where the Schur function $P_d(t_1,\cdots,t_d)$ is defined as:
$$P_d(t_1,\cdots,t_d)=\sum_{m=1}^d\sum_{d_1+\cdots + d_m=d(d_i>0)}\frac {t_{d_1}\cdots t_{d_m}} {m!}  \qquad  (P_0=1)  $$
or equivalently
$$exp\left (\sum_{k=1}^{\infty}t_{k}z^k\right )= \sum_{k=0}^{\infty}P_k\left (t_{1},\cdots, t_{k}\right )z^k
 \eqno {(1.4)}$$

It was also proved in [4, Theorem 6.10] that $Tr_k(\mathbb {T})$ is the sum of $k$-th powers of all the eigenvalues of $\mathbb {T}$. In the following Lemma 1.1 and Theorem 1.1, we show that this important result can also be simply proved by using the Schur function as follows.

\vskip 0.48cm

\noindent {\bf Lemma 1.1 }: Let $a_0=1$, and $\sum_{i=0}^na_iz^{n-i}=(z-\lambda_1)\cdots (z-\lambda_n)$ be a monic polynomial of degree $n$ with $n$ roots $\lambda_1,\cdots ,\lambda_n$.
Let $t_k=-\frac {1} {k}(\sum_{j=1}^n \lambda_j^k)$, and $g(z)=\sum_{k=1}^{\infty}t_kz^k$. Then we have
$e^{g(z)}=\sum_{k=0}^na_kz^k$.

\vskip 0.28cm

\noindent {\bf Proof}: Differentiating both sides of the equation $g(z)=\sum_{k=1}^{\infty}t_kz^k$, we have (since $\sum_{k=1}^{\infty}t_kz^k$ has a positive radius of convergence)
$$g'(z)=-\sum_{j=1}^n \lambda_j\left ( \sum_{k=1}^{\infty}(\lambda_jz)^{k-1}  \right )=\left (\sum_{j=1}^n \ln (1-\lambda_jz) \right )'$$
Also $g(0)=0$, so we have $g(z)=\sum_{j=1}^n \ln (1-\lambda_jz) $, and thus
$$e^{g(z)}=\prod_{j=1}^n(1-\lambda_jz)=z^n\left ( \sum_{i=0}^na_iz^{-(n-i)} \right )=\sum_{k=0}^na_kz^k$$
\qed

\vskip 0.28cm

\noindent {\bf Theorem 1.1 ([4])}: Let $\mathbb {T}$ be a tensor of order $m\ge 2$ and dimension $n$, let $Tr_k(\mathbb {T})$ be defined as in (1.3). Let $\lambda_1,\cdots ,\lambda_d$ be all the eigenvalues of $\mathbb {T}$ (where $d=n(m-1)^{n-1}$). Then we have:
$$Tr_k(\mathbb {T})=\sum_{j=1}^d \lambda_j^k  \eqno {(1.5)}$$

\vskip 0.28cm

\noindent {\bf Proof}: Let $\phi_{\mathbb {T}}(z)=\sum_{i=0}^da_iz^{d-i}=(z-\lambda_1)\cdots (z-\lambda_d)$ be the characteristic polynomial of the tensor $\mathbb {T}$. Let  $t_k=-\frac {1} {k}(\sum_{j=1}^d \lambda_j^k)$, and $g(z)=\sum_{k=1}^{\infty}t_kz^k$. Then by Lemma 1.1 we have $e^{g(z)}=\sum_{k=0}^da_kz^k$.
\vskip 0.18cm

On the other hand, by  [3,4,8] we know that the coefficient $a_k$ of the characteristic polynomial of the tensor $\mathbb {T}$ is $a_k=P_k(-\frac {Tr_1(\mathbb {T})} {1},\cdots, -\frac {Tr_k(\mathbb {T})} {k})$, and
$P_k(-\frac {Tr_1(\mathbb {T})} {1},\cdots, -\frac {Tr_k(\mathbb {T})} {k})=0$ when $k>d$. Thus by (1.4) we also have
$$exp\left (\sum_{k=1}^{\infty}-\frac {Tr_k(\mathbb {T})} {k}z^k\right )= \sum_{k=0}^{\infty}P_k\left (-\frac {Tr_1(\mathbb {T})} {1},\cdots, -\frac {Tr_k(\mathbb {T})} {k}\right )z^k=\sum_{k=0}^da_kz^k$$
Comparing this with the expression for $e^{g(z)}$, we obtain that
$$\sum_{k=1}^{\infty}-\frac {Tr_k(\mathbb {T})} {k}z^k=g(z)=\sum_{k=1}^{\infty}t_kz^k$$
and thus $-\frac {Tr_k(\mathbb {T})} {k}=t_k=-\frac {1} {k}(\sum_{j=1}^d \lambda_j^k)$. From this (1.5) follows.
\qed

\vskip 0.28cm

Notice that the current formula (1.3) for the trace $Tr_k(\mathbb {T})$ (given in [3] and [8]) involves the differential operators and auxiliary matrix $A$, so it is quite difficult and complex to use it to study the traces. And it is hoped that some more explicit formulas for the trace $Tr_k(\mathbb {T})$ can be obtained (as was mentioned in the final remarks of [4]).

\vskip 0.28cm

In this paper, we will give in \S 2 and \S 4 some new formulas for the trace $Tr_k(\mathbb {T})$ in terms of some graph theoretical parameters. These formulas do not involve the differential operators and auxiliary matrix. In \S 3, we give three applications of the trace formula given in \S 2 in the study of the spectra of uniform hypergraphs. Firstly we give a characterization (in terms of the traces of tensors) of the $k$-uniform hypergraphs whose spectra are $k$-symmetric, thus give an answer to a question raised in [3]. Secondly we generalize the $k$-partite and hm-bipartite hypergraphs to p-hm bipartite hypergraphs, and prove that the spectra of this class of hypergraphs are $k$-symmetric if $p$ and $k$ are coprime. This result is a common generalization of the results [3, Theorem 4.2] and [5, Proposition 3.1]. Thirdly, we answer a question raised in [5] about the relation between the Laplacian spectrum and signless Laplacian spectrum of a $k$-uniform hypergraph. In \S 5, we use the new trace formulas given in \S 4 to give a simplified proof of a formula for $Tr_2(\mathbb {T})$, and discuss the possible expression for $Tr_3(\mathbb {T})$.

\vskip 1.88cm

\section{A new graph theoretical formula for the trace $Tr_k(\mathbb {T})$ }

\vskip 0.28cm

In this section, we first use the weighted associated digraph $D(A)$ of a matrix $A$ to give a graph theoretical expression for $tr(A^r)$ in Lemma 2.2. Then we derive a graph theoretical formula (2.11) for the trace $Tr_k(\mathbb {T})$ of a tensor $\mathbb {T}$ in Theorem 2.1.

\vskip 0.28cm

A multi-set is a collection of elements which allows the repeated elements. In this paper, if a multi set $A$ contains $s$ distinct elements $a_1,\cdots,a_s$ with the multiplicities $r_1,\cdots,r_s$ respectively, then we write
$$A=a_1^{r_1}\cdots a_s^{r_s}$$

\vskip 0.28cm

\noindent {\bf Lemma 2.1}: Let $a_1,\cdots ,a_n$ and $b_1,\cdots ,b_n$ be nonnegative integers with $a_1+\cdots +a_n=b_1+\cdots +b_n$. Then we have:

$$\frac {\partial^{a_1+\cdots +a_n}} {(\partial x_1)^{a_1}(\partial x_2)^{a_2}\cdots (\partial x_n)^{a_n}} (x_1^{b_1}x_2^{b_2}\cdots x_n^{b_n})=\left\{
                            \begin{array}{cc}
                              b_1!\cdots b_n! &   \mbox{if } a_i=b_i \ (i\in \{1,\cdots ,n\})  \\
                              0 &  \mbox{otherwise}\\
                            \end{array}
                          \right. \eqno {(2.1)}$$

\noindent {\bf Proof}. If some $a_i\ne b_i$, then by the condition $a_1+\cdots +a_n=b_1+\cdots +b_n$ we must have some $a_j>b_j$. Thus in this case the left side of (2.1) is zero. The case when $a_i=b_i$ for all $i=1,\cdots ,n$ is obvious.   \qed

\vskip 0.28cm
In the following, we write $[n]=\{1,\cdots ,n\}$.

\vskip 0.28cm

\noindent {\bf Definition 2.1}: Let $A=(a_{ij})$ be a matrix of order $n$. Then the weighted associated digraph $D(A)$ of $A$ is the digraph with vertex set $V=[n]$ such that there is an arc $(i,j)$ in $D(A)$ if and only if $a_{ij}\ne 0$, and in this case, the arc $(i,j)$ has a weight $a_{ij}$. The weight of a walk $W$ in $D(A)$, denoted by $a(W)$, is the product of the weights of all arcs of $W$ (here the arc set of $W$ is regarded as a multi set). Let $\mathbf {W}_r(D(A))$ be the set of all closed walks of length $r$ in $D(A)$.

\vskip 0.18cm

The following lemma gives a graph theoretical formula for the trace $tr(A^r)$ in terms of the weighted associated digraph $D(A)$ of a matrix $A$.

\vskip 0.28cm

\noindent {\bf Lemma 2.2}: Let $A=(a_{ij})$ be a matrix of order $n$ with the weighted associated digraph $D(A)$. Then we have:
$$tr(A^r)=\sum_{W\in \mathbf {W}_r(D(A))}a(W)  \eqno {(2.2)}$$

\noindent {\bf Proof}. We have

$$tr(A^r)= \sum_{i_1,\cdots ,i_r=1}^na_{i_1i_2}a_{i_2i_3}\cdots  a_{i_ri_1} =\sum_{W\in \mathbf {W}_r(D(A))}a(W)  $$
\qed

\vskip 0.18cm

For a tensor $\mathbb {H}=(h_{i_1i_2\cdots i_m})$ of order $m$ and dimension $n$, we write:

$$h_{i_1i_2\cdots i_m}=h_{i_1\alpha } \qquad  (\mbox {where} \ \alpha =i_2\cdots i_m \in [n]^{m-1})$$

\noindent For an integer $d>0$, we define:

$$\mathcal{F}_d=\{\left ((i_1,\alpha_1 ),\cdots , (i_d,\alpha_d ) \right ) \ | \ 1\le i_1\le \cdots \le i_d\le n; \  \  \alpha_1 ,\cdots ,\alpha_d \in [n]^{m-1} \}  \eqno {(2.3)}$$

\noindent For nonnegative integers $d_1,\cdots , d_n$ with $d_1+\cdots + d_n=d>0$, we also define:

$$\mathcal{F}_{d_1,\cdots , d_n}=\{((i_1,\alpha_1 ),\cdots , (i_d,\alpha_d ))\in \mathcal{F}_{d} \ | \ \{i_1, \cdots , i_d\}= 1^{d_1}\cdots n^{d_n} \}  \eqno {(2.4)}$$
Then we obviously have
$$\mathcal{F}_d = \bigcup_{d_1+\cdots + d_n=d} \mathcal{F}_{d_1,\cdots , d_n} $$

To prove our formulas for $Tr_d(\mathbb {T})$, we need the following elementary formula on the transformation of reversing order between sums and products.

\vskip 0.28cm

\noindent {\bf Lemma 2.3}:
For a tensor $\mathbb {H}=(h_{i_1i_2\cdots i_m})$ of order $m$ and dimension $n$ with
$h_{i_1i_2\cdots i_m}$ denoted by $h(i_1, \alpha)$,   $(\alpha =i_2\cdots i_m \in [n]^{m-1})$, we have:

$$\prod_{i=1}^n\left ( \sum_{y_i\in [n]^{m-1}}h(i,y_i)\right )^{d_i}= \sum_{((i_1,\alpha_1 ),\cdots , (i_d,\alpha_d ))\in \mathcal{F}_{d_1,\cdots , d_n}}\prod_{j=1}^d h(i_j, \alpha_j )   \eqno {(2.5)}$$

\noindent {\bf Proof}. We have
$$
\begin{aligned}
\displaystyle &  \prod_{i=1}^n\left (\sum_{y_i\in [n]^{m-1}}h(i,y_i)\right )^{d_i}     \\
&=  \sum_{y_{ij}\in [n]^{m-1} \ (i=1,\cdots,n, \ j=1,\cdots,d_i)} h(1,y_{11})\cdots h(1,y_{1d_1}) h(2,y_{21})\cdots h(2,y_{2d_2}) \cdots h(n,y_{n1})\cdots h(n,y_{nd_n})  \\
&= \sum_{((i_1,\alpha_1 ),\cdots , (i_d,\alpha_d ))\in \mathcal{F}_{d_1,\cdots , d_n}}\prod_{j=1}^d h(i_j, \alpha_j )    \\
\end{aligned}  $$
\qed

\vskip 0.18cm

Now for $F=((i_1,\alpha_1 ),\cdots , (i_d,\alpha_d ))\in \mathcal{F}_d$ and the tensor $\mathbb {H}=(h_{i_1i_2\cdots i_m})$ of order $m$ and dimension $n$, write
$$\pi_F(\mathbb {H})=\prod_{j=1}^d h(i_j, \alpha_j ) \eqno {(2.6)}$$
Also for a given $F\in \mathcal{F}_d$, there exist unique nonnegative integers $d_1,\cdots , d_n$ with  $d_1+\cdots + d_n=d$ such that $F\in \mathcal{F}_{d_1,\cdots , d_n}$. In this case, we write (for any one variable function $g(x)$):
$$g(F)=\prod_{i=1}^n g(d_i)  \qquad (F\in \mathcal{F}_{d_1,\cdots , d_n}) \eqno {(2.7)}$$

\vskip 0.18cm

Using these notations and Lemma 2.3, we further have

\vskip 0.28cm

\noindent {\bf Lemma 2.4}:
$$\sum_{d_1+\cdots + d_n=d}\prod_{i=1}^n\left (g(d_i)\left (\sum_{y_i\in [n]^{m-1}}h(i,y_i)\right )^{d_i}\right )=\sum_{F\in \mathcal{F}_d}g(F)\pi_F(\mathbb {H}) \eqno {(2.8)}$$

\noindent {\bf Proof}. From (2.5), (2.6) and (2.7) we have

$$
\begin{aligned}
\displaystyle &  \sum_{d_1+\cdots + d_n=d}\prod_{i=1}^n\left (g(d_i)\left (\sum_{y_i\in [n]^{m-1}}h(i,y_i)\right )^{d_i}\right ) \\
& =\sum_{d_1+\cdots + d_n=d}\left (\prod_{i=1}^ng(d_i)\right )\left ( \sum_{((i_1,\alpha_1 ),\cdots , (i_d,\alpha_d ))\in \mathcal{F}_{d_1,\cdots , d_n}}\prod_{j=1}^d h(i_j, \alpha_j ) \right )   \\
& =\sum_{d_1+\cdots + d_n=d}\left (\prod_{i=1}^ng(d_i)\right )\sum_{F\in \mathcal{F}_{d_1,\cdots , d_n}}\pi_F(\mathbb {H}) = \sum_{d_1+\cdots + d_n=d}\sum_{F\in \mathcal{F}_{d_1,\cdots , d_n}}\left (\prod_{i=1}^ng(d_i)\right )\pi_F(\mathbb {H})  \\
&= \sum_{F\in \mathcal{F}_d}g(F)\pi_F(\mathbb {H})
\\
\end{aligned}  $$
\qed

Now we introduce some more graph theoretical notations. We first assume that all the digraphs considered here have the vertex set $V=[n]$, and may have loops (arcs of the form $(i,i)$) and multiple arcs (such digraphs are called ``multi-digraphs"). Thus there may be several arcs from vertex $i$ to vertex $j$. For convenience, we use $D_n$ to denote the complete digraph of order $n$ with the arc set $E(D_n)=[n]\times [n]$.

\vskip 0.18cm

For an arc multi-set $E$, we use $V(E)$ to denote the set of vertices incident to some arc of $E$.

\vskip 0.18cm

In the following, when we mention an arc multi-set $E$, we always mean that $E$ is a multi-set each of whose element is in $[n]\times [n]$. Namely, $V(E)\subseteq [n]$. Also, for each vertex $i\in [n]$, let $d_E^+(i)$ and $d_E^-(i)$ be the outdegree and indegree of $i$ in the arc set $E$, respectively.

\vskip 0.18cm

A (multi) digraph is called a "balanced digraph", if the outdegree and indegree of each vertex are equal.
It is not difficult to see that, if $W$ is a closed walk (of some digraph with vertex set $V=[n]$), then $W$ (as a multi-digraph) is a "balanced digraph".

\vskip 0.38cm

\noindent {\bf Definition 2.2}: Let $E$ be an arc multi-set (with $V(E)\subseteq [n]$). Then:

\vskip 0.18cm
\noindent (1) Let $b(E)$ be the product of the factorials of the multiplicities of all the arcs of $E$.

\vskip 0.18cm
\noindent (2) Let $c(E)$ be the product of the factorials of the outdegrees of all the vertices in the arc set $E$.

\vskip 0.18cm
\noindent (3) Let $\mathbf{W}(E)$ be the set of all (directed) closed walks $W$ with the arc multi-set $E(W)=E$.

\vskip 0.38cm

\noindent {\bf Definition 2.3}: Let $F=((i_1,\alpha_1 ),\cdots , (i_d,\alpha_d ))\in \mathcal{F}_d$, where $(i_j, \alpha_j )\in [n]^m \ (j=1,\cdots, d)$. Then

\vskip 0.18cm
\noindent (1) Let $E(F)=\bigcup _{j=1}^dE_j(F)$ (in the sense of the union of multi-sets), where $E_j(F)$ is the arc multi-set as following:
$$E_j(F)=\{ (i_j,v_1), (i_j,v_2), \cdots (i_j,v_{m-1}) \} \qquad \mbox {if} \  \ \alpha_j=(v_1,\cdots,v_{m-1})  $$
Thus $E(F)$ is also an arc multi-set.

\vskip 0.18cm
\noindent (2) Let $b(F)=b(E(F))$ be the product of the factorials of the multiplicities of all the arcs of $E(F)$.

\vskip 0.18cm
\noindent (3)  Let $c(F)=c(E(F))$ be the product of the factorials of the outdegrees of all the vertices in the arc set $E(F)$.

It is easy to see that if $F\in \mathcal{F}_{d_1,\cdots , d_n}$, then $d_{E(F)}^+(i)=d_i(m-1)$. Thus in this case we have $c(F)=\prod_{i=1}^n(d_i(m-1))!$.

\vskip 0.18cm
\noindent (4) Let $\mathbf{W}(F)=\mathbf{W}(E(F))$ be the set of all closed walks $W$ with the arc multi-set $E(W)=E(F)$. It is obvious that the length of each walk $W$ in $\mathbf{W}(F)$ is $|E(W)|=|E(F)|=d(m-1)$.

\vskip 0.18cm
\noindent (5) Let the differential operator $\partial (F)=\prod_{j=1}^d \frac {\partial} {\partial a_{i_j\alpha_j}}$, where

$$\frac {\partial} {\partial a_{i\alpha}}=\prod_{k=1}^{m-1} \frac {\partial} {\partial a_{is_k}} \qquad (\mbox {if}  \  \  \alpha = (s_1,\cdots,s_{m-1}) \in [n]^{m-1} )$$
Here $a_{ij} \ (i,j=1,\cdots,n)$ are viewed as distinct independent variables.

\vskip 0.18cm

Now let $\mathbb {T}=(t_{i_1i_2\cdots i_m})$ be a tensor of order $m$ and dimension $n$, where
$t_{i_1i_2\cdots i_m}=t_{i_1\alpha }$   $ \ (\alpha =i_2\cdots i_m \in [n]^{m-1})$.
Take the tensor $\mathbb {H}$ in (2.8) as $h_{i\alpha}=t_{i\alpha}\frac {\partial} {\partial a_{i\alpha}}$ (viewed as an element in some operator algebra). Then $\pi_F(\mathbb {H})=\pi_F(\mathbb {T})\partial (F)$, and from (2.8) we have:

$$\sum_{d_1+\cdots + d_n=d}\prod_{i=1}^n \frac {1} {(d_i(m-1))!}\left (\sum_{y_i\in [n]^{m-1}}t_{iy_i}\frac {\partial} {\partial a_{iy_i}}\right )^{d_i}=\sum_{F\in \mathcal{F}_d}\frac {1} {c(F)}\pi_F(\mathbb {T})\partial (F) \eqno {(2.9)}$$

We also have the following formula about the action of the differential operator $\partial (F)$ on $tr(A^{d(m-1)})$.

\vskip 0.28cm

\noindent {\bf Lemma 2.5}: For $F=((i_1,\alpha_1 ),\cdots , (i_d,\alpha_d ))\in \mathcal{F}_d$, let
$\partial (F)=\prod_{j=1}^d \frac {\partial} {\partial a_{i_j\alpha_j}}$ be defined as in Definition 2.2, and
$A=(a_{ij})$ be a matrix of order $n$, where $a_{ij} \ (i,j=1,\cdots,n)$ are distinct independent variables.  Then we have

$$\partial (F)(tr(A^{d(m-1)}))=b(F)|\mathbf{W}(F)| \qquad  (F\in \mathcal{F}_d)  \eqno {(2.10)}$$

\noindent {\bf Proof}. By using the trace formula (2.2) for matrix $A^r$, we have

$$\partial (F) tr(A^{d(m-1)})=\sum_{W\in \mathbf {W}_{d(m-1)}(D(A))}\partial (F)a(W)  $$

Now for fixed $W\in \mathbf {W}_{d(m-1)}(D(A))$ and $F\in \mathcal{F}_d$, we know by Lemma 2.1 that $\partial (F)a(W)\ne 0$ if and only if the arc multi-sets $E(W)=E(F)$, namely $W\in \mathbf{W}(F)$, and in this case $\partial (F)a(W)=b(F)$ by Lemma 2.1. Thus we have

$$\partial (F) tr(A^{d(m-1)})=\sum_{W\in \mathbf{W}(F)}\partial (F)a(W) =\sum_{W\in \mathbf{W}(F)}b(F)=b(F)|\mathbf{W}(F)|  $$
\qed

Now we are ready to prove our first graph theoretical trace formula.

\vskip 0.28cm

\noindent {\bf Theorem 2.1}: Let $\mathbb {T}=(t_{i_1i_2\cdots i_m})$ be a tensor of order $m$ and dimension $n$. Then we have
$$Tr_d(\mathbb {T})= (m-1)^{n-1} \sum_{F\in \mathcal{F}_d}\frac {b(F)} {c(F)}\pi_F(\mathbb {T})|\mathbf{W}(F)|  \eqno {(2.11)}$$
where (the graph theoretical parameters) $b(F)$ $c(F)$ and $|\mathbf{W}(F)|$ only depend on the arc set $E(F)$, and are independent of the tensor $\mathbb {T}$.

\vskip 0.18cm

\noindent {\bf Proof}. By (2.9) we have

$$\sum_{d_1+\cdots + d_n=d}\prod_{i=1}^n \frac {1} {(d_i(m-1))!}\left (\sum_{y_i\in [n]^{m-1}}t_{iy_i}\frac {\partial} {\partial a_{iy_i}}\right )^{d_i}tr(A^{d(m-1)})=\sum_{F\in \mathcal{F}_d}\frac {1} {c(F)}\pi_F(\mathbb {T})  \partial (F) tr(A^{d(m-1)})$$
Substituting (2.10) into the above equation, we obtain:
$$\sum_{d_1+\cdots + d_n=d}\prod_{i=1}^n \frac {1} {(d_i(m-1))!}\left (\sum_{y_i\in [n]^{m-1}}t_{iy_i}\frac {\partial} {\partial a_{iy_i}}\right )^{d_i}tr(A^{d(m-1)})=\sum_{F\in \mathcal{F}_d}\frac {b(F)} {c(F)}\pi_F(\mathbb {T})|\mathbf{W}(F)|  \eqno {(2.12)}$$
Multiplying both sides of (2.12) by $(m-1)^{n-1}$, we obtain our trace formula (2.11).
\qed

\vskip 1.28cm

\section{ Some applications in the study of spectra of hypergraphs}

In this section, we give some applications of the trace and our trace formula (2.11) in the study of the spectra (and Laplacian spectra) of hypergraphs. Firstly we give a characterization (in terms of the traces of tensors) of the $k$-uniform hypergraphs whose spectra are $k$-symmetric. Secondly we generalize the $k$-partite and hm-bipartite hypergraphs to p-hm bipartite hypergraphs, and prove that the spectra of this class of hypergraphs are $k$-symmetric if $p$ and $k$ are coprime. Finally, we answer a question raised in [5] about the relation between the Laplacian spectrum and signless Laplacian spectrum of a $k$-uniform hypergraph.

\vskip 0.18cm

\vskip 0.28cm

A hypergraph $H=(V,E)$ is called $k$-uniform if every edge of $H$ contains exactly $k$ vertices. The adjacency tensor of $H$ (under certain ordering of vertices) is the order $k$ dimension $n$ tensor $\mathbb {A}=\mathbb {A}_H$ with the following entries ([3]):
$$a_{i_1i_2\cdots i_k}=\left \{\begin{array}{cc}
                               \frac {1} {(k-1)!} &   \mbox{if }  \ \{i_1,i_2,\cdots ,i_k\} \in E(H) \\
                               0 &   \mbox{otherwise}\\
                            \end{array}
                          \right.$$
The characteristic polynomial and spectrum of a uniform hypergraph $H$ are that of its adjacency tensor $\mathbb {A}$.

\vskip 0.2cm

Let $\mathbb {D}=\mathbb {D}_H$ be the degree diagonal tensor of $H$ (its $i$-th diagonal element is the degree of the vertex $i$), then the tensor $\mathbb {L}=\mathbb {D}-\mathbb {A}$ is called the Laplacian tensor of $H$, and $\mathbb {Q}=\mathbb {D}+\mathbb {A}$ is called the signless Laplacian tensor of $H$. The Laplacian spectrum and signless Laplacian spectrum of $H$ are defined to be the spectrum of $\mathbb {L}$ and $\mathbb {Q}$, respectively.

\vskip 0.18cm

The spectrum of a tensor or a $k$-uniform hypergraph is said to be $k$-symmetric, if this spectrum is invariant under a rotation of an angle $2\pi/k$ in the complex plane.

\vskip 0.18cm

In [3, Theorem 4.2], Cooper and Dutle proved that the spectrum of a $k$-partite $k$-uniform hypergraph is $k$-symmetric. They also proposed a problem in [3] which asks to characterize those hypergraphs whose spectra are $k$-symmetric. In the following Theorem 3.1, we will give a characterization of the $k$-uniform hypergraphs whose spectra are $k$-symmetric in terms of the traces of its adjacency tensor, and then we will give an application of this result (together with the new trace formula (2.11)) in Theorem 3.2 to show that the spectra of the class of p-hm hypergraphs are $k$-symmetric when $p,k$ are coprime.

\vskip 0.28cm

\noindent {\bf Theorem 3.1.} Let $H$ be a $k$-uniform hypergraph, $\mathbb {A}=\mathbb {A}_H$ be its adjacency tensor, and $\phi_{\mathbb {A}}(\lambda )=\sum_{j=0}^ra_j\lambda ^{r-j} \ (r=n(k-1)^{n-1})$ be the characteristic polynomial of $\mathbb {A}$ and $H$. Then the following three conditions are equivalent:

\vskip 0.18cm
\noindent (1). The spectrum of $\mathbb {A}$ (and $H$) is $k$-symmetric.

\vskip 0.18cm
\noindent (2). If $d$ is not a multiple of $k$, then the coefficient $a_d$ (of the codegree $d$ term in the characteristic polynomial $\phi _{\mathbb {A}} (\lambda )$) is zero. Namely,
 there exist some integer $t$ and some polynomial $f$, such that $\phi _{\mathbb {A}} (\lambda )=\lambda ^t f(\lambda ^k)$.

\vskip 0.18cm
\noindent (3). If $d$ is not a multiple of $k$, then $Tr_d(\mathbb {A})=0$.

\vskip 0.18cm

\noindent {\bf Proof}. We will show that $(1)\Longleftrightarrow (2)\Longleftrightarrow (3)$.

\vskip 0.18cm
$(2)\Longrightarrow (1)$: This is obvious from the expression $\phi _{\mathbb {A}} (\lambda )=\lambda ^t f(\lambda ^k)$.

\vskip 0.18cm

$(1)\Longrightarrow (2)$:
Let $\varepsilon =e^{2\pi i/k}$ be the $k$-th primitive root of unity. Then (1) implies that $\phi _{\mathbb {A}} (\varepsilon \lambda )= \varepsilon ^r\phi _{\mathbb {A}} (\lambda )$. From this we have
$$\sum_{d=0}^ra_d\varepsilon ^{r-d}\lambda ^{r-d} =\sum_{d=0}^ra_d\varepsilon ^{r}\lambda ^{r-d}$$
Thus we have $a_d\varepsilon ^{r-d}= a_d\varepsilon ^{r}$, or $a_d(\varepsilon ^{d}-1)=0$.

Now if $d$ is not a multiple of $k$, then $\varepsilon ^{d}-1\ne 0$. So in this case we have $a_d=0$.

\vskip 0.18cm
$(2)\Longrightarrow (3)$: By (2), we may write $\phi _{\mathbb {A}} (\lambda )$ as:
$$\phi _{\mathbb {A}} (\lambda )=\lambda ^t(\lambda ^k-c_1^k)\cdots (\lambda ^k-c_s^k)$$
Let $P=(p_{ij})$ be the circulant permutation matrix of order $k$ (where $p_{ij}=1$ if and only if $j\equiv i+1 (mod \ k)$). If $d$ is not a multiple of $k$, then all the diagonal entries of $P^d$ are zero, thus $tr(P^d)=0$.

Also, we have $\phi _{cP} (\lambda )= \lambda ^k-c^k$. So if $\mu_1,\cdots, \mu_k$ are the $k$ roots of $\lambda ^k-c^k=0$ (i.e., the $k$ eigenvalues of the matrix $cP$), then $\mu_1^d,\cdots, \mu_k^d$ will be the $k$ eigenvalues of the matrix $(cP)^d$. Thus we have
$$\mu_1^d+\cdots+ \mu_k^d=tr((cP)^d)$$
Therefore by Theorem 1.1 we have
$$Tr_d(\mathbb {A})=\sum_{j=1}^r \lambda _j^d = tr((c_1P)^d)+\cdots +tr((c_sP)^d)=0 \qquad \mbox {(if $d$ is not a multiple of $k$)}$$

\vskip 0.18cm
$(3)\Longrightarrow (2)$:  By [3,4,8], we have
$$a_d=P_d(-\frac {Tr_1(\mathbb {T})} {1},\cdots, -\frac {Tr_d(\mathbb {T})} {d}) \eqno {(3.1)}$$
where $P_d(t_1,\cdots,t_d)$ is the Schur function defined as in Section 1 ([3,8]):
$$P_d(t_1,\cdots,t_d)= \sum_{m=1}^d\sum_{d_1+\cdots+d_m=d(d_i>0)}\frac {t_{d_1}\cdots t_{d_m}} {m!}   \eqno {(3.2)}$$

Now suppose that $a_d\ne 0$. Then by (3.1) and (3.2) we see that there exist some positive integers $d_1,\cdots,d_m$ with $d_1+\cdots+d_m=d$ such that
$$Tr_{d_1}(\mathbb {A})\cdots Tr_{d_m}(\mathbb {A})\ne 0  \eqno {(3.3)}$$

By condition (3), we see that (3.3) implies that $d_1,\cdots,d_m$ are all multiples of $k$. Thus $d=d_1+\cdots+d_m$ is also a multiple of $k$,  this proves (2).
\qed

\vskip 0.28cm

In [5], Hu and Qi defined the ($k$-uniform) hm-bipartite hypergraphs (which is a generalization of the $k$-partite hypergraphs studied in [3]), and  proved that a number $\lambda_0$ is an eigenvalue of an hm-bipartite hypergraph $H$ if and only if $\lambda_0 e^{2\pi i/k}$ is an eigenvalue of $H$. In the following, we further generalize hm-bipartite hypergraphs to $p$-hm bipartite hypergraphs, and prove the $k$-symmetry of the spectra of $p$-hm bipartite hypergraphs when $p,k$ are coprime, thus generalize both the results [3, Theorem 4.2] and [5, Proposition 3.1].

\vskip 0.28cm

\noindent {\bf Definition 3.1 }. Let $H=(V,E)$ be a nontrivial $k$-uniform hypergraph. It is called $p$-hm bipartite if $V$ can be partitioned into $V=V_1\cup V_2$, where $V_1$ and $V_2$ are nonempty and disjoint, such that every edge of $H$ intersects $V_1$ with exactly $p$ vertices.

\vskip 0.2cm

The hm-bipartite hypergraphs defined in [5] is a special case $p=1$ of $p$-hm bipartite hypergraphs. Also, the cored hypergraphs defined in [6] (every edge contains a vertex of degree one) is a special class of $1$-hm bipartite hypergraphs.

\vskip 0.2cm

In order to prove the $k$-symmetry of the spectra of $p$-hm bipartite hypergraphs when $p,k$ are coprime, we need the following Lemma 3.1 which is in some sense an equivalent version of [3, Theorem 3.12].

\vskip 0.38cm

Let $F=\left ((i_1,\alpha_1 ),\cdots , (i_d,\alpha_d ) \right )\in \mathcal{F}_d$ (where each component of $F$ is an element of $[n]^m$), and let $i\in [n]$. Let $d_i(F)$ be the number of times that the index $i$ appears in $F$ as the primary index (i.e., the first index in some component of $F$), and $q_i(F)$ be the number of times that the index $i$ appears in $F$ as the non-primary index. Let $p_i(F)=d_i(F)+q_i(F)$ be the total number of times that the index $i$ appears in $F$. Then it was defined in [3] that $F$ is called $m$-valent, if for each $i\in [n]$, $p_i(F)$ is a multiple of $m$.

\vskip 0.28cm

\noindent {\bf Lemma 3.1 }. Let $F=\left ((i_1,\alpha_1 ),\cdots , (i_d,\alpha_d ) \right )\in \mathcal{F}_d$ (where each component of $F$ is an element of $[n]^m$). If $\mathbf{W}(F)\ne \phi$, then $F$ is $m$-valent.

\vskip 0.18cm

\noindent {\bf Proof}. Take $W\in \mathbf{W}(F)$ to be a closed walk with $E(W)=E(F)$. Then we have $d_W^+(i)=d_W^-(i)$ for each vertex $i\in [n]$ (since $E(W)$ is balanced). Now by the definition of $E(F)$, we can see that $d_W^+(i)=(m-1)d_i(F)$ and $d_W^-(i)=q_i(F)$. Thus we have $q_i(F)=(m-1)d_i(F)$, and so $p_i(F)=d_i(F)+q_i(F)=d_i(F)+(m-1)d_i(F)=md_i(F)$, which is a multiple of $m$.
\qed

\vskip 0.28cm

Now let
$$\mathcal{F}'_d=\{F\in \mathcal{F}_d \ | \ F \ \mbox {is $m$-valent} \ \}  \eqno {(3.4)}$$
Then from Lemma 3.1 we can see that (2.11) can be rewritten as (in terms of $\mathcal{F}'_d$):
$$Tr_d(\mathbb {T})= (m-1)^{n-1} \sum_{F\in \mathcal{F}'_d}\frac {b(F)} {c(F)}\pi_F(\mathbb {T})|\mathbf{W}(F)|  \eqno {(3.5)}$$
since for those $F\in \mathcal{F}_d \backslash \mathcal{F}'_d$, we have $|\mathbf{W}(F)|=0$ by Lemma 3.1.

\vskip 0.18cm

Now we apply (3.5) to prove the following theorem.

\vskip 0.28cm

\noindent {\bf Theorem 3.2 }. Let $H=(V,E)$ be a nontrivial $k$-uniform $p$-hm bipartite hypergraph with $p,k$ coprime, then the spectrum of $H$ is $k$-symmetric.

\vskip 0.18cm

\noindent {\bf Proof}. Let $\mathbb {A}=\mathbb {A}_H$ be the adjacency tensor of $H$. By Theorem 3.1, we only need to show that $\mathbb {A}$ satisfies the condition (3) of Theorem 3.1.

Let $V_1$ and $V_2$ be as in Definition 3.1.

\vskip 0.18cm
Suppose that $Tr_d(\mathbb {A})\ne 0$ for some positive integer $d$. Then by the formula (3.5), there exists some $F\in \mathcal{F}'_d$ such that $\pi_F(\mathbb {A})\ne 0$. Thus the $d$ components of $F$ corresponds to $d$ edges $\{e_1,\cdots, e_d\}$ (with repetition allowed) of the hypergraph $H$, and that $F$ is $k$-valent by (3.4).

Let $E_0=\{e_1,\cdots, e_d\}$, and for each vertex $v$ of $H$, let $d_{E_0}(v)$ be the degree of $v$ in the sub-hypergraph of $H$ induced by the edge subset $E_0$. Then by the $k$-valent property of $F$, we see that all $d_{E_0}(v)$ are multiples of $k$.

On the other hand, by Definition 3.1 (and the same idea as in [3]), we have $pd=\sum_{v\in V_1}d_{E_0}(v)$, which is a multiple of $k$. So $d$ is also a multiple of $k$ since $p,k$ are coprime. This proves that $\mathbb {A}$ satisfies the condition (3) of Theorem 3.1.
\qed

\vskip 0.18cm

Now we consider an application of the formula (2.11) in the study of the Laplacian spectra of hypergraphs. In [5], it was asked whether the Laplacian spectrum and signless  Laplacian spectrum equal or not for a $k$-uniform hypergraph with odd $k\ge 3$. By using our trace formula (2.11), we are able to answer this question in the following theorem.

\vskip 0.18cm

\vskip 0.18cm

\noindent {\bf Theorem 3.3}: Let $H$ be a nontrivial $k$-uniform hypergraph with odd $k\ge 3$. Then its Laplacian spectrum and signless  Laplacian spectrum do not equal.

\vskip 0.18cm

\noindent {\bf Proof}. Let $\mathbb {A}=\mathbb {A}_H$ be the adjacency tensor, and $\mathbb {D}$ be the diagonal degree tensor of $H$. Let $\mathbb {L}=\mathbb {D}-\mathbb {A}$ be the Laplacian tensor of $H$. Then we obviously have $\mathbb {D}+\mathbb {A}=|\mathbb {L}|$, where $|\mathbb {L}|$ is obtained from $\mathbb {L}$ by taking the absolute values entrywisely. Now by Theorem 2.1, we have:
$$Tr_k(\mathbb {D}-\mathbb {A})= Tr_k(\mathbb {L})= (k-1)^{n-1} \sum_{F\in \mathcal{F}_k}\frac {b(F)} {c(F)}\pi_F(\mathbb {L})|\mathbf{W}(F)|  \eqno {(3.6)}$$
and
$$Tr_k(\mathbb {D}+\mathbb {A})= Tr_k(|\mathbb {L}|)= (k-1)^{n-1} \sum_{F\in \mathcal{F}_k}\frac {b(F)} {c(F)}\pi_F(|\mathbb {L}|)|\mathbf{W}(F)| = (k-1)^{n-1} \sum_{F\in \mathcal{F}_k}\left | \frac {b(F)} {c(F)}\pi_F(\mathbb {L})|\mathbf{W}(F)| \right |  \eqno {(3.7)}$$
Namely, each term of the right side of (3.7) is the absolute value of the corresponding term  of the right side of (3.6).

\vskip 0.18cm

Now we want to show that there exists some $F_0\in \mathcal{F}_k$ such that $\pi_{F_0}(\mathbb {L})<0$ and
$|\mathbf{W}(F_0)|>0$. For this purpose, take an edge $e=\{i_1,\cdots,i_k\}$ of the hypergraph $H$ and take
$$F_0=\left ( ( i_1,i_2\cdots,i_k), (i_2,i_3,\cdots, i_k,i_1),\cdots, (i_k,i_1,\cdots, i_{k-1})    \right )\in \mathcal{F}_k $$
Then we have $\pi_{F_0}(\mathbb {L})=\pi_{F_0}(-\mathbb {A}) =(-1)^k\frac {1} {((k-1)!)^k}<0$ since $k$ is odd.

\vskip 0.18cm

On the other hand, we have $E(F_0)=\{(i,j) \ | \ i,j\in \{i_1,i_2\cdots,i_k\} \ \mbox {and} \ i\ne j\}$, which means that the digraph (on $n$ vertices) induced by the arc set $E(F_0)$ is isomorphic to the complete digraph $D_k$ (without loops) which is both strongly connected and balanced ($d_{D_k}^+(v)=d_{D_k}^-(v)=k-1$ for all vertices $v$ in $D_k$). Thus by the criterion for the "directed Eulerian graphs" ([1]), we conclude that there exists a directed closed walk $W\in \mathbf{W}(F_0)$ with $E(W)=E(F_0)$. So we have $|\mathbf{W}(F_0)|>0$.

Now for this $F_0$, we have $\frac {b(F_0)} {c(F_0)}\pi_{F_0}(\mathbb {L})|\mathbf{W}(F_0)| <0$, since $b(F)>0$ and $c(F)>0$ for all $F\in \mathcal{F}_k$. From this we see that at least one term in the summation of the righthand side of (3.6) is negative, so the righthand sides of (3.6) and (3.7) do not equal. Consequently, we have $Tr_k(\mathbb {D}-\mathbb {A})\ne Tr_k(\mathbb {D}+\mathbb {A})$. From this and Theorem 1.1 we conclude that $\mathbb {D}-\mathbb {A}$ and $\mathbb {D}+\mathbb {A}$ have the different spectra.
\qed

\vskip 1.28cm

\section{Some other formulas for the trace $Tr_k(\mathbb {T})$ }

\vskip 0.28cm

In this section, we give some more formulas for the trace $Tr_k(\mathbb {T})$, and consider some examples and applications. In order to obtain these trace formulas, we need to introduce some more graph theoretical notations.

\vskip 0.38cm

\noindent {\bf Definition 4.1}: Let $n,d,r$ be fixed positive integers. Let $\mathbf {E}_{d,r}(n)$ be the set of arc multi-sets $E$ (with $V(E)\subseteq [n]$) satisfying the following three conditions:

\vskip 0.18cm
\noindent (1) $|E|=dr$ (in the sense of multi set).

\vskip 0.18cm
\noindent (2) The arc multi-set $E$ is balanced (i.e., for each vertex $i\in [n]$, the outdegree $d_E^+(i)$ and the indegree $d_E^-(i)$ are equal).

\vskip 0.18cm
\noindent (3) The outdegree $d_E^+(i)$ of every vertex $i\in [n]$ is a multiple of $r$.

\vskip 0.18cm

By the condition (3) of the above definition, we see that if $E\in \mathbf {E}_{d,r}(n)$, then $d_E^+(i)\ge r$ if $d_E^+(i)>0$. On the other hand, recall that $V(E)=\{i\in [n] \ | \ d_E^+(i)>0\}$. So $dr=|E|=\sum_{i=1}^nd_E^+(i)\ge |V(E)|r$. Thus we have

$$E\in \mathbf {E}_{d,r}(n)\Longrightarrow |V(E)|\le d  \eqno {(4.1)} $$
where equality holds if and only if every vertex $i$ with $d_E^+(i)>0$ has $d_E^+(i)=d_E^-(i)=r$.

\vskip 0.38cm

\noindent {\bf Lemma 4.1}: Let $n,m,d$ be fixed positive integers and $F\in \mathcal{F}_d$. If $\mathbf{W}(F)\ne \phi$, then we have $E(F)\in \mathbf {E}_{d,m-1}(n)$.

\vskip 0.28cm

\noindent {\bf Proof}. From the definitions of $\mathcal{F}_d$ and $E(F)$, it is easy to see that $F\in \mathcal{F}_d$ implies $|E(F)|=d(m-1)$.

Also, by the hypothesis we have $\mathbf{W}(E(F))=\mathbf{W}(F)\ne \phi$, which means that there is a closed walk $W\in \mathbf{W}(E(F))$ with $E(F)$ as its arc multi set. Thus $E(F)$ is balanced, since $E(W)$ is.

Furthermore, $F\in \mathcal{F}_d$ implies that $F\in \mathcal{F}_{d_1,\cdots , d_n}$ for some nonnegative integers $d_1,\cdots , d_n$ with $d_1+\cdots +d_n=d$. Then it is easy to see that $d_{E(F)}^+(i)=d_i(m-1)$ which is a multiple of $(m-1)$ for all $i\in [n]$. Thus $E(F)$ satisfies all the three conditions in Definition 4.1 (with $r=m-1$), so $E(F)\in \mathbf {E}_{d,m-1}(n)$.
\qed

\vskip 0.38cm

By Lemma 4.1, we have:
$$\{F\in \mathcal{F}_d \  | \ \mathbf{W}(F)\ne \phi \}\subseteq \bigcup_{E\in \mathbf {E}_{d,m-1}(n)} \{F\in \mathcal{F}_d \  | \ E(F)=E \} \  \subseteq  \  \mathcal{F}_d  \eqno {(4.2)} $$
Thus the trace formula (2.11) can be further written as:
$$
\begin{aligned}
\displaystyle  Tr_d(\mathbb {T}) & = (m-1)^{n-1} \sum_{F\in \mathcal{F}_d}\frac {b(F)} {c(F)}\pi_F(\mathbb {T})|\mathbf{W}(F)|   \\
& = (m-1)^{n-1} \sum_{E(F)\in \mathbf {E}_{d,m-1}(n)}\frac {b(F)} {c(F)}\pi_F(\mathbb {T})|\mathbf{W}(F)|   \\
&=  (m-1)^{n-1} \sum_{E\in \mathbf {E}_{d,m-1}(n)}\sum_{F\in \mathcal{F}_d , E(F)=E}\frac {b(E)} {c(E)}\pi_F(\mathbb {T})|\mathbf{W}(E)|
\\
\end{aligned} \eqno {(4.3)} $$
We also have:

\vskip 0.38cm

\noindent {\bf Lemma 4.2}: Let $b(E)$ and $c(E)$ be defined as in Definition 2.2. Then for any arc multi-set

\noindent $E= \bigcup_{i=1}^n \bigcup_{j=1}^n (i,j)^{r_{ij}} \in \mathbf {E}_{d,m-1}(n)$, where the multiplicity of the arc $(i,j)$ in $E$ is $r_{ij}$, and $d_E^+(i)=\sum_{j=1}^nr_{ij}$, we have

$$|\{F\in \mathcal{F}_d \  | \ E(F)=E \}|=\frac {c(E)} {b(E)} \eqno {(4.4)} $$

\vskip 0.28cm

\noindent {\bf Proof}. For each fixed $i$, we list all the elements of $E$ with the initial vertex $i$ as following:
$$(i,1),\cdots,(i,1); \cdots ; (i,n),\cdots,(i,n) \qquad \mbox {where there are $r_{ij}$ many $(i,j)$'s } \ (j=1,\cdots,n) \eqno {(4.5)} $$
By using the formula for the number of permutations with repetition, we know that the number of the permutations of the elements in (4.5) is the following multi-binomial coefficient:

$$\left (
\begin{array}{c}
r_{i1}+\cdots +r_{in}\\
r_{i1},\cdots ,r_{in}\\
\end{array} \right )=\frac {(r_{i1}+\cdots +r_{in})!} {r_{i1}!\cdots r_{in}!}  $$
Thus we have

$$|\{F\in \mathcal{F}_d \  | \ E(F)=E \}|=\prod_{i=1}^n \frac {(r_{i1}+\cdots +r_{in})!} {r_{i1}!\cdots r_{in}!}
= \frac {\prod_{i=1}^n(d_E^+(i))!} {\prod_{i=1}^n\prod_{j=1}^nr_{ij}!} =\frac {c(E)} {b(E)} $$
\qed

\vskip 0.28cm

Now for each $E\in \mathbf {E}_{d,m-1}(n)$, let
$$\pi _E(\mathbb{T})= \sum_{F\in \mathcal{F}_d , E(F)=E}\pi_F(\mathbb {T})\eqno {(4.6)}$$
and
$$\overline{\pi _E(\mathbb{T})}=\frac {\sum_{F\in \mathcal{F}_d , E(F)=E}\pi_F(\mathbb {T})} {|\{F\in \mathcal{F}_d \  | \ E(F)=E \}|}=\frac {b(E)} {c(E)}\pi _E(\mathbb{T})  \eqno {(4.7)}$$
Thus $\overline{\pi _E(\mathbb{T})}$ is the average value of all those values $\pi_F(\mathbb {T})$ with $F\in \mathcal{F}_d$ and $E(F)=E$. Using this and the equation (4.3), we can now obtain the following two trace formulas.

\vskip 0.28cm

\noindent {\bf Theorem 3.1}: Let $\mathbb {T}$ be a tensor of order $m$ and dimension $n$. Then we have:
$$Tr_d(\mathbb {T})= (m-1)^{n-1} \sum_{E\in \mathbf {E}_{d,m-1}(n)}\frac {b(E)} {c(E)}\pi _E(\mathbb{T})|\mathbf{W}(E)|  \eqno {(4.8)}$$
and
$$Tr_d(\mathbb {T})= (m-1)^{n-1} \sum_{E\in \mathbf {E}_{d,m-1}(n)}\overline{\pi _E(\mathbb{T})}|\mathbf{W}(E)|  \eqno {(4.9)}$$

\vskip 0.28cm

\noindent {\bf Proof}. (4.8) follows directly from (4.3) and (4.6), while (4.9) follows directly from (4.8) and (4.7).
\qed

\vskip 0.28cm

Now we look at some examples.

\vskip 0.18cm

\noindent {\bf Example 4.1}. Let the tensor $\mathbb {T}=A=(a_{ij})$ be a matrix of order $n$ (i.e., $m=2$). Then we have $Tr_d(\mathbb {T})= tr(A^d)$.

\vskip 0.28cm

\noindent {\bf Proof}. First we have $m-1=1$. Next, for any arc multi-set $E\in \mathbf {E}_{d,1}(n)$, let
$$a(E)=\prod_{e\in E}a(e)  \qquad  \mbox {where} \  \ a(e)=a_{ij} \  \mbox {if} \ e=(i,j)$$
be the weight of the arc multi-set $E$ in the weighted associated digraph $D(A)$.

Now for any $F=((i_1,j_1),\cdots, (i_d,j_d)\in \mathcal{F}_d$ with $E(F)=E$ and $W\in \mathbf {W}(E)$, we have
$\pi_F(\mathbb {T})=a_{i_1j_1}\cdots a_{i_dj_d}=a(E(F))=a(E)=a(W)$, where $a(W)$ is defined in Definition 2.1. Thus we have $\overline{\pi _E(\mathbb{T})}=a(E)$, and so
$$\sum_{W\in \mathbf {W}(E)}a(W)=a(E)\sum_{W\in \mathbf {W}(E)}1= \overline{\pi _E(\mathbb{T})}|\mathbf{W}(E)|  \eqno {(4.10)}$$

Also, we have
$$\mathbf {W}_d(D(A))  \subseteq \bigcup_{E\in \mathbf {E}_{d,1}(n)}\mathbf {W}(E)$$
and if $W\in \bigcup_{E\in \mathbf {E}_{d,1}(n)}\mathbf {W}(E)\backslash \mathbf {W}_d(D(A))$, then $a(W)=0$. Thus by (2.2), (4.9) and (4.10) we have

$$
\begin{aligned}
\displaystyle   tr(A^d) &= \sum_{W\in \mathbf {W}_d(D(A))}a(W)\\
& = \sum_{E\in \mathbf {E}_{d,1}(n)}\sum_{W\in \mathbf {W}(E)}a(W)\\
& = \sum_{E\in \mathbf {E}_{d,1}(n)}\overline{\pi _E(\mathbb{T})}|\mathbf{W}(E)| \\
& =Tr_d(\mathbb {T})\\
\end{aligned}  $$
\qed

\vskip 0.38cm

\noindent {\bf Example 4.2}. Let $\mathbb {J}$ be the tensor of order $m$ and dimension $n$ with all elements 1. Let

$$\mathbf {W}_{d,m-1}(n)=\{W  \ \mbox {is a closed walk} \ | \ E(W)\in \mathbf {E}_{d,m-1}(n)\}$$
Then we have $Tr_d(\mathbb {J})= (m-1)^{n-1} |\mathbf {W}_{d,m-1}(n)|$, where
$$|\mathbf {W}_{d,m-1}(n)|=\sum_{d_1+\cdots+d_n=d}\frac {(d(m-1))!} {\prod_{i=1}^n(d_i(m-1))!}$$

\vskip 0.28cm

\noindent {\bf Proof}. By Definition 2.2 we know that $W\in \mathbf {W}(E)$ if and only if $E(W)=E$. So by the definition of $\mathbf {W}_{d,m-1}(n)$ we have
$$\mathbf {W}_{d,m-1}(n)=\bigcup_{E\in \mathbf {E}_{d,m-1}(n)}\mathbf {W}(E) \ , \qquad \mbox {so} \quad  \
|\mathbf {W}_{d,m-1}(n)|=\sum_{E\in \mathbf {E}_{d,m-1}(n)}|\mathbf {W}(E)|.  $$

Since all the elements of $\mathbb {J}$ are 1, we have $\overline{\pi _E(\mathbb{J})}=1$ for any $E\in \mathbf {E}_{d,m-1}(n)$. Thus by (4.9) and the above equation we have

$$Tr_d(\mathbb {J})= (m-1)^{n-1} \sum_{E\in \mathbf {E}_{d,m-1}(n)}|\mathbf{W}(E)|=(m-1)^{n-1}|\mathbf {W}_{d,m-1}(n)|$$
By using the formula for the permutations with repetitions, and viewing a closed walk as a sequence of vertices, we can see that
$$|\mathbf {W}_{d,m-1}(n)|=\sum_{d_1+\cdots+d_n=d}\frac {(d(m-1))!} {\prod_{i=1}^n(d_i(m-1))!}$$

\noindent ({\bf Note}: Here $|\mathbf {W}_{d,m-1}(n)|$ is a combinatorial parameter which only depends on $n,m$ and $d$, and is independent of the tensors.)    \qed

\vskip 1.28cm

\section{The expressions of $Tr_2(\mathbb {T})$ and $Tr_3(\mathbb {T})$}

In this section we show how our trace formulas can be used in the study of $Tr_2(\mathbb {T})$ and $Tr_3(\mathbb {T})$. First, we use the formula (4.8) and (4.9) to give a simplified proof of a formula of $Tr_2(\mathbb {T})$ in [4].

\vskip 0.28cm

\noindent {\bf Theorem 5.1 ([4])}. Let $\mathbb {T}$ be a tensor with order $m$ and dimension $n$. Then we have

$$Tr_2(\mathbb {T})=(m-1)^{n-1}\left [\sum_{i=1}^nt_{ii\cdots i}^2 +\sum_{i<j}\sum_{s=1}^{m-1}\frac {2s} {m-1}    \left ( \sum_{\{i_2,\cdots i_m\}=j^si^{m-1-s}}t_{ii_2\cdots i_m}   \right ) \left ( \sum_{\{j_2,\cdots j_m\}=i^sj^{m-1-s}}t_{jj_2\cdots j_m}   \right )           \right ]  $$

\vskip 0.28cm

\noindent {\bf Proof}. We use the formula (4.9).

For each $E\in \mathbf {E}_{2,m-1}(n)$, let
$V(E)= \{i\in [n] \ | \ d_E^+(i)>0\}$ as before. Then $|V(E)|\le 2$ by (4.1), so we can divide the set $\mathbf {E}_{2,m-1}(n)$ into two subsets as $\mathbf {E}_{2,m-1}(n)=\mathbf {E}_1\bigcup \mathbf {E}_2$, where $E\in \mathbf {E}_k$ if and only if $|V(E)|=k$ (for $k=1,2$).

Thus we can further write
$$\mathbf {E}_1=\{E(1),\cdots, E(n)\}, \qquad  \mbox {and}  \qquad  \mathbf {E}_2=\bigcup _{i<j} \mathbf {E}(i,j) \eqno {(5.1)} $$
where $V(E(i))=\{i\}$ (thus $E(i)=(i,i)^{2(m-1)}$), and $V(E)=\{i,j\}$ for each $E\in \mathbf {E}(i,j)$. Furthermore, for each $1\le i<j\le n$, we can write
$$\mathbf {E}(i,j)=\{E_0(i,j),E_1(i,j),\cdots, E_{m-1}(i,j)\}$$
where (as arc multi set)
$$E_s(i,j)=(i,j)^s(j,i)^s(i,i)^{m-1-s}(j,j)^{m-1-s}  \quad (0\le s\le m-1, \ i<j)  \eqno {(5.2)} $$

Now for $E=E(i)=(i,i)^{2(m-1)}$, we have $|\mathbf{W}(E)|=1$ and $\overline{\pi _E(\mathbb{T})}=t_{i\cdots i}^2$. So
$$\sum_{E\in \mathbf {E}_1}\overline{\pi _E(\mathbb{T})}|\mathbf{W}(E)| =\sum_{i=1}^nt_{i\cdots i}^2 \eqno {(5.3)} $$

For $E=E_s(i,j)$, we can verify from (5.2) that
$$b(E)=(s!(m-1-s)!)^2,\qquad \qquad c(E)=((m-1)!)^2$$
and
$$\pi_E(\mathbb {T})=\sum_{F\in \mathcal{F}_2 , E(F)=E_s(i,j)}\pi_F(\mathbb {T})=\left ( \sum_{\{i_2,\cdots i_m\}=j^si^{m-1-s}}t_{ii_2\cdots i_m}   \right ) \left ( \sum_{\{j_2,\cdots j_m\}=i^sj^{m-1-s}}t_{jj_2\cdots j_m}   \right ) \eqno {(5.4)} $$

Now we consider $\mathbf{W}(E)$ for $E=E_s(i,j)$. If $W\in \mathbf{W}(E)$, then the initial vertex of $W$ is either $i$ or $j$. If the initial vertex of $W$ is $i$, then there are $\left (
\begin{array}{c}
m-1\\
s\\
\end{array} \right )$  different orderings of the $m-1$ arcs in $W$ starting from $i$, since there are $s$ arcs $(i,j)$ among these $m-1$ arcs. On the other hand, among the $m-1$ arcs in $W$ starting from $j$, the last arc must be $(j,i)$ since the terminal vertex of $W$ is also $i$. Thus there are only $\left (
\begin{array}{c}
m-2\\
s-1\\
\end{array} \right )$ different orderings of the remaining $m-2$ arcs in $W$ starting from $j$. The similar arguments apply if the initial vertex of $W$ is $j$. Therefore we have $|\mathbf{W}(E)|= 2\left (
\begin{array}{c}
m-1\\
s\\
\end{array} \right )\left (
\begin{array}{c}
m-2\\
s-1\\
\end{array} \right )$
for $E=E_s(i,j)$.

Combining this with the expressions for $b(E)$ and $c(E)$, we have for $E=E_s(i,j)$ that
$$\frac {b(E)} {c(E)} |\mathbf{W}(E)|=\frac {(s!(m-1-s)!)^2} {((m-1)!)^2}2\left (
\begin{array}{c}
m-1\\
s\\
\end{array} \right )\left (
\begin{array}{c}
m-2\\
s-1\\
\end{array} \right )=\frac {2\left (
\begin{array}{c}
m-1\\
s\\
\end{array} \right )\left (
\begin{array}{c}
m-2\\
s-1\\
\end{array} \right )} {\left (
\begin{array}{c}
m-1\\
s\\
\end{array} \right )^2}=\frac {2s} {m-1}  \eqno {(5.5)}$$

Finally, by using the formula (4.9) together with (5.3), (5.4), (5.5) and (4.7), we have

$$
\begin{aligned}
\displaystyle  Tr_2(\mathbb {T}) & =(m-1)^{n-1} \sum_{E\in \mathbf {E}_{2,m-1}(n)}\overline{\pi _E(\mathbb{T})}|\mathbf{W}(E)| \\
& = (m-1)^{n-1}\left ( \sum_{E\in \mathbf {E}_1}\overline{\pi _E(\mathbb{T})}|\mathbf{W}(E)|+\sum_{E\in \mathbf {E}_2}\overline{\pi _E(\mathbb{T})}|\mathbf{W}(E)| \right ) \\
&= (m-1)^{n-1}\left ( \sum_{E\in \mathbf {E}_1}\overline{\pi _E(\mathbb{T})}|\mathbf{W}(E)|+\sum_{i<j}\sum_{s=0}^{m-1}\overline{\pi _{E_s(i,j)}(\mathbb{T})}|\mathbf{W}(E_s(i,j))| \right )
\\
&= (m-1)^{n-1}\left ( \sum_{i=1}^nt_{i\cdots i}^2+\sum_{i<j}\sum_{s=0}^{m-1}\frac {b(E_s(i,j))} {c(E_s(i,j))}|\mathbf{W}(E_s(i,j))|\pi_{E_s(i,j)}(\mathbb {T})   \right )  \\
&=  (m-1)^{n-1}\left [\sum_{i=1}^nt_{ii\cdots i}^2 +\sum_{i<j}\sum_{s=1}^{m-1}\frac {2s} {m-1}    \left ( \sum_{\{i_2,\cdots i_m\}=j^si^{m-1-s}}t_{ii_2\cdots i_m}   \right ) \left ( \sum_{\{j_2,\cdots j_m\}=i^sj^{m-1-s}}t_{jj_2\cdots j_m}   \right )           \right ]  \\
\end{aligned}  $$
Here the first equality follows from (4.9), the fourth equality follows from (5.3) and (4.7), and the last equality follows from (5.4) and (5.5).
\qed

\vskip 0.68cm

Finally, we consider $Tr_3(\mathbb {T})$. Similarly as for $Tr_2(\mathbb {T})$,  we can divide the set $\mathbf {E}_{3,m-1}(n)$ into three subsets as $\mathbf {E}_{3,m-1}(n)=\mathbf {E}_1\bigcup \mathbf {E}_2\bigcup \mathbf {E}_3$, where $E\in \mathbf {E}_k$ if and only if $|V(E)|=k$ (for $k=1,2,3$), since $|V(E)|\le 3$ for each $E\in \mathbf {E}_{3,m-1}(n)$. Then we consider the following three cases.

\vskip 0.28cm

\noindent {\bf Case 1 (For $|V(E)|=1$)}. We obviously have $\mathbf {E}_1=\{E(1),\cdots, E(n)\}$,
where $V(E(i))=\{i\}$, and so  $E(i)=(i,i)^{3(m-1)}$. Thus we have $|\mathbf{W}(E(i))|=1$ and $\overline{\pi _E(\mathbb{T})}=t_{i\cdots i}^3$. Therefore
$$\sum_{E\in \mathbf {E}_1}\overline{\pi _E(\mathbb{T})}|\mathbf{W}(E)| =\sum_{i=1}^nt_{i\cdots i}^3 \eqno {(5.6)} $$

\vskip 0.28cm

\noindent {\bf Case 2 (For $|V(E)|=2$)}. We have $\mathbf {E}_2=\bigcup _{i\ne j} \mathbf {E}(i,j)$, where $V(E)=\{i,j\}$ for each $E\in \mathbf {E}(i,j)$. Furthermore, we have
$$\mathbf {E}(i,j)=\{E_0(i,j),E_1(i,j),\cdots, E_{m-1}(i,j)\}$$
where
$$E_s(i,j)=(i,j)^s(j,i)^s(i,i)^{2(m-1)-s}(j,j)^{m-1-s}  \quad (0\le s\le m-1, \ i\ne j)  \eqno {(5.7)} $$

Now for $E=E_s(i,j)$, we have
$$b(E)=(s!)^2(m-1-s)!(2(m-1)-s)!,\qquad \qquad c(E)=(m-1)!(2(m-1))!$$
and similarly as in the case for $Tr_2(\mathbb {T})$, we also have $|\mathbf{W}(E)|= \left (
\begin{array}{c}
m-1\\
s\\
\end{array} \right )\left (
\begin{array}{c}
2(m-1)-1\\
s-1\\
\end{array} \right )+\left (
\begin{array}{c}
2(m-1)\\
s\\
\end{array} \right )\left (
\begin{array}{c}
m-2\\
s-1\\
\end{array} \right )$
for $E=E_s(i,j)$ (where the first term corresponds to those closed walks with initial vertex $j$, and the second term corresponds to those closed walks with initial vertex $i$).

Combining this with the expressions for $b(E)$ and $c(E)$, we have for $E=E_s(i,j)$ that

$$
\begin{aligned}
\displaystyle  & \frac {b(E)} {c(E)} |\mathbf{W}(E)|= \\
& =\frac {1} {\left (
\begin{array}{c}
m-1\\
s\\
\end{array} \right )\left (
\begin{array}{c}
2(m-1)\\
s\\
\end{array} \right )} \left ( \left (
\begin{array}{c}
m-1\\
s\\
\end{array} \right )\left (
\begin{array}{c}
2(m-1)-1\\
s-1\\
\end{array} \right )+\left (
\begin{array}{c}
2(m-1)\\
s\\
\end{array} \right )\left (
\begin{array}{c}
m-2\\
s-1\\
\end{array} \right ) \right ) \\
& =\frac {s} {2(m-1)}+\frac {s} {m-1}=\frac {3s} {2(m-1)} \\
\end{aligned}  \eqno {(5.8)} $$

For $E=E_s(i,j)$, we also have

$$\begin{aligned}
\displaystyle   \pi_E(\mathbb {T}) & =\sum_{F\in \mathcal{F}_3 , E(F)=E_s(i,j)}\pi_F(\mathbb {T}) \\
& =\left ( \sum_{\{i_2,\cdots i_m,k_2,\cdots,k_m\}=j^si^{2(m-1)-s}}t_{ii_2\cdots i_m} t_{ik_2\cdots k_m}  \right ) \left ( \sum_{\{j_2,\cdots j_m\}=i^sj^{m-1-s}}t_{jj_2\cdots j_m}   \right )\\
\end{aligned}  $$
Thus we have

$$\begin{aligned}
\displaystyle  & \sum_{E\in \mathbf {E}_2}\overline{\pi _E(\mathbb{T})}|\mathbf{W}(E)|= \\
 & =\sum_{i\ne j}\sum_{s=0}^{m-1}\frac {3s} {2(m-1)}\left ( \sum_{\{i_2,\cdots i_m,k_2,\cdots,k_m\}=j^si^{2(m-1)-s}}t_{ii_2\cdots i_m} t_{ik_2\cdots k_m}  \right ) \left ( \sum_{\{j_2,\cdots j_m\}=i^sj^{m-1-s}}t_{jj_2\cdots j_m}   \right )\\
\end{aligned}   \eqno {(5.9)} $$

\vskip 0.28cm

\noindent {\bf Case 3 (For $|V(E)|=3$)}. We have $\mathbf {E}_3=\bigcup _{i< j<k} \mathbf {E}(i,j,k)$, where $V(E)=\{i,j,k\}$ for each $E\in \mathbf {E}(i,j,k)$.

Now for fixed $1\le i< j<k \le n$ and each $E\in \mathbf {E}(i,j,k)$, we have $d_E^+(i)=d_E^+(j)=d_E^+(k)=d_E^-(i)=d_E^-(j)=d_E^-(k)=m-1$. Let $p,q,r,s$ be the multiplicities of the arcs $(i,j)$, $(j,k)$, $(k,i)$ and $(j,i)$ in $E$, respectively. Then $E$ must have the following form:
$$
\begin{aligned}
\displaystyle  & E=E(i,j,k;p,q,r,s) \\
&:=(i,j)^p(j,k)^q(k,i)^r(j,i)^s(i,k)^{r+s-p}(k,j)^{q+s-p}(i,i)^{m-1-s-r}
(j,j)^{m-1-s-q}(k,k)^{m-1+p-r-s-q} \\
\end{aligned}  \eqno {(5.10)} $$
Thus we have
$$\mathbf {E}(i,j,k)=\{E(i,j,k;p,q,r,s) \ | \ 0\le p,q,r,s \le m-1, \  \mbox {and all the multiplicities in (5.10) are nonnegative} \} \eqno {(5.11)} $$

Now for $E=E(i,j,k;p,q,r,s)$, we have $c(E)=((m-1)!)^3$, and
$$b(E)=p!q!r!s!(r+s-p)!(q+s-p)!(m-1-s-r)!(m-1-s-q)!(m-1+p-r-s-q)!$$
So
$$\frac {b(E)} {c(E)}  =\frac {1} {\left (
\begin{array}{c}
m-1\\
s, r,  m-1-s-r\\
\end{array} \right )\left (
\begin{array}{c}
m-1\\
p,  q+s-p, m-1-s-q\\
\end{array} \right )\left (
\begin{array}{c}
m-1\\
q,  r+s-p,  m-1+p-r-s-q\\
\end{array} \right )} \eqno {(5.12)} $$

Now let $w(m;p,q,r,s)$ be the number of closed walks $W$ with the arc multi set $E(W)=E(i,j,k;p,q,r,s)$. Then
$w(m;p,q,r,s)$ is a purely graph theoretical parameter (which only involves the digraphs with three vertices), and is independent of $n,i,j,k$, and independent of the tensor $\mathbb {T}$. Then we have
$$|W(E)|=w(m;p,q,r,s) \qquad (\mbox {if} \  \  E=E(i,j,k;p,q,r,s)) \eqno {(5.13)} $$
(for those values $0\le p,q,r,s \le m-1$ such that some multiplicities in (5.10) are negative, we would have
$|W(E)|=w(m;p,q,r,s)=0$.)
\vskip 0.28cm

Also, let
$$t(i,j,k;p,q,r,s)=\left ( \sum_{\{i_2,\cdots i_m\}=i^*j^pk^{r+s-p}}t_{ii_2\cdots i_m}   \right ) \left ( \sum_{\{j_2,\cdots j_m\}=i^sj^*k^q}t_{jj_2\cdots j_m}   \right )\left ( \sum_{\{k_2,\cdots k_m\}=i^rj^{q+s-p}k^*}t_{kk_2\cdots k_m}   \right ) \eqno {(5.14)} $$
(where the *'s mean suitable numbers such that the total multiplicities of the multi sets $\{i_2,\cdots i_m\}$,  $\{j_2,\cdots j_m\}$ and $\{k_2,\cdots k_m\}$ are all $m-1$.)
Then we have
$$\pi _{E(i,j,k;p,q,r,s)}(\mathbb{T})=t(i,j,k;p,q,r,s) \eqno {(5.15)} $$

From (5.11), (5.12), (5.13) and (5.15) we have
$$\sum_{E\in \mathbf {E}_3}\frac {b(E)} {c(E)}\pi _E(\mathbb{T})|\mathbf{W}(E)| =\sum_{i<j<k}\sum_{p=0}^{m-1}\sum_{q=0}^{m-1}\sum_{r=0}^{m-1}\sum_{s=0}^{m-1}\frac {w(m;p,q,r,s)t(i,j,k;p,q,r,s)} {\left (
\begin{array}{c}
m-1\\
s, \ r, \ *\\
\end{array} \right )\left (
\begin{array}{c}
m-1\\
p, \ q+s-p, \ *\\
\end{array} \right )\left (
\begin{array}{c}
m-1\\
q, \ r+s-p, \ *\\
\end{array} \right )} \eqno {(5.16)} $$
(where the *'s mean suitable numbers which make the sums of the corresponding three numbers equal to $m-1$.)

\vskip 0.28cm

Combining (5.6), (5.9) and (5.16), we finally have:

\vskip 0.28cm

\noindent {\bf Theorem 5.2}. Let $\mathbb {T}$ be the tensor with order $m$ and dimension $n$. Then we have

$$
\begin{aligned}
\displaystyle  &\frac {Tr_3(\mathbb {T})} {(m-1)^{n-1}}= \\
&\sum_{i=1}^nt_{ii\cdots i}^3 + \sum_{i\ne j}\sum_{s=0}^{m-1}\frac {3s} {2(m-1)}\left ( \sum_{\{i_2,\cdots i_m,k_2,\cdots,k_m\}=j^si^{2(m-1)-s}}t_{ii_2\cdots i_m} t_{ik_2\cdots k_m}  \right ) \left ( \sum_{\{j_2,\cdots j_m\}=i^sj^{m-1-s}}t_{jj_2\cdots j_m}   \right ) \\
&+ \sum_{i<j<k}\sum_{p=0}^{m-1}\sum_{q=0}^{m-1}\sum_{r=0}^{m-1}\sum_{s=0}^{m-1}\frac {w(m;p,q,r,s)t(i,j,k;p,q,r,s)} {\left (
\begin{array}{c}
m-1\\
s, \ r, \ *\\
\end{array} \right )\left (
\begin{array}{c}
m-1\\
p, \ q+s-p, \ *\\
\end{array} \right )\left (
\begin{array}{c}
m-1\\
q, \ r+s-p, \ *\\
\end{array} \right )} \\
\end{aligned}     \eqno {(5.17)}     $$
where $w(m;p,q,r,s)$ and $t(i,j,k;p,q,r,s)$ are defined in (5.13) and (5.14).

\vskip 0.28cm

\noindent {\bf Proof}. By (4.9) and $\mathbf {E}_{3,m-1}(n)=\mathbf {E}_1\bigcup \mathbf {E}_2\bigcup \mathbf {E}_3$ we have:
$$\frac {Tr_3(\mathbb {T})} {(m-1)^{n-1}}= \sum_{E\in \mathbf {E}_1}\overline{\pi _E(\mathbb{T})}|\mathbf{W}(E)| + \sum_{E\in \mathbf {E}_2}\overline{\pi _E(\mathbb{T})}|\mathbf{W}(E)|+\sum_{E\in \mathbf {E}_3}\overline{\pi _E(\mathbb{T})}|\mathbf{W}(E)|$$
Substituting (5.6), (5.9) and (5.16) into the above equation, we obtain (5.17).
\qed

\vskip 0.28cm

\noindent {\bf Remark 5.1}.

\vskip 0.18cm
\noindent (1) For those values $0\le p,q,r,s \le m-1$ such that some multiplicities in (5.10) are negative, we have $w(m;p,q,r,s)=0$. So adding or deleting some terms corresponding to these values $p,q,r,s$ (or those values $p,q,r,s$ such that the corresponding $w(m;p,q,r,s)=0$)  will not affect the value of the total sum in the expression of our formula.

\vskip 0.18cm
\noindent (2) When $m$ is small, the graph theoretical parameter $w(m;p,q,r,s)$ can be computed directly or by using a computer.

\vskip 0.18cm
\noindent (3) A further comment on the term $\sum_{E\in \mathbf {E}_3}\frac {b(E)} {c(E)}\pi _E(\mathbb{T})|\mathbf{W}(E)|$ in (5.16).

Theoretically, if we obtain a formula of $\sum_{E\in \mathbf {E}_3}\frac {b(E)} {c(E)}\pi _E(\mathbb{T})|\mathbf{W}(E)|$ for the tensor $\mathbb{T}$ of order $m$ and dimension $n=3$, then we can obtain a formula of $\sum_{E\in \mathbf {E}_3}\frac {b(E)} {c(E)}\pi _E(\mathbb{T})|\mathbf{W}(E)|$ for a general tensor $\mathbb{T}$ of order $m$ and dimension $n$ by doing the following two steps:

\vskip 0.18cm
\noindent {\bf Step 1}: In each term of the formula of $\sum_{E\in \mathbf {E}_3}\frac {b(E)} {c(E)}\pi _E(\mathbb{T})|\mathbf{W}(E)|$ (for the tensor $\mathbb{T}$ of order $m$ and dimension $n=3$), replace all the subscripts 1 by $i$, 2 by $j$, and 3 by $k$, and keep the coefficients of that term unchanged.

\vskip 0.18cm
\noindent {\bf Step 2}: Add $\sum_{i<j<k}$ at the beginning of that part (for $E\in \mathbf {E}_3$) of the formula.

\vskip 0.8cm
\noindent{\bf References}

\bibliographystyle{acm}

\end{document}